\documentclass{article}

\title{Doing Algebra over an Associative Algebra}

\author{Nathan BeDell \\ nbedell@tulane.edu}

\usepackage[a4paper, total={6in, 8in}]{geometry}
\usepackage{graphicx}
\usepackage{amsmath}
\usepackage{amssymb}
\usepackage{mathrsfs}
\usepackage{amsfonts}
\usepackage[utf8]{inputenc}
\usepackage[english]{babel}
\usepackage{amsthm}
\usepackage{amscd}
\usepackage{tikz}
\usepackage{bbm}

\newcommand{\Acal}{\mathcal{A}}

\newcommand{\Hcal}{\mathcal{H}}
\newcommand{\Pcal}{\mathcal{P}}
\newcommand{\zd}{\mathbf{zd}}

\newcommand{\Ann}{\mathrm{Ann}}

\theoremstyle{definition}
\newtheorem{definition}{Definition}[section]

\newtheorem{theorem}[definition]{Theorem}
\newtheorem{example}[definition]{Example}
\newtheorem{lemma}[definition]{Lemma}
\newtheorem{corollary}[definition]{Corollary}
\newtheorem{proposition}[definition]{Proposition}

\theoremstyle{remark}

\theoremstyle{remark}

\theoremstyle{remark}

\theoremstyle{remark}

\begin{document}
	
	\maketitle
	
	\begin{abstract}
		A finite-dimensional unital and associative algebra over $\mathbb{R}$, or what we shall call simply ``an algebra'' in this paper for short, generalities the construction by which we derive the complex numbers by ``adjoining an element $i$'' to $\mathbb{R}$ and imposing the relation $i^2 = -1$. In this paper, we examine some of the elementary algebraic properties of such algebras, how they break-down when compared to standard grade-school algebra, and discuss how such properties are relevant to other areas of our research regarding algebras, such as the $\Acal$-calculus and the theory of $\Acal$-ODEs. 
	\end{abstract}
	
\section{Preface}

Ever since Gerolamo Cardano stumbled upon the idea of the complex numbers in the 16th century, mathematicians have been thinking about alternate systems of numbers which generalize the reals. In this paper, we look at different methods of constructing such number systems, which are technically known as ``finite-dimensional unital associative algebras'', but in this paper we will simply refer to as ``algebras''. (Not to be confused with the subject matter of ``algebra'', which is the sense the word is first used in the title of this paper.)

Of course, since the field (pun intended) of ring theory is already well-developed in modern mathematics, one might wonder what new might come of a paper such as this one. And the answer is: perhaps not much, at least in terms of technical results. The novelty of this paper is that it gives a more hands-on, elementary approach (when possible) to the study of such algebraic systems, more akin to what is called ``Algebra I \& II'' in the United States, except in the context of algebras. In other words: how do the basics of algebra (i.e. divisibility, factorization theory) hold up in a more general context; that of finite-dimensional associative unital algebras (i.e. ``algebras'')? Hence, ``Doing Algebra over an Associative Algebra''. Of course, if one delves too deeply into these questions, one quickly finds themselves back in modern algebra, so we merely scratch the surface here.  

Thus, the main focus of this paper will be to attempt to inculcate an elementary understanding of how the principles of ``usual grade-school algebra'' break down in the presence of new features, such as zero-divisors, and especially nilpotent elements, and how, without using (too much) high-powered abstract theorems from modern algebra, one might understand and characterize this break-down. 

We will more-or-less assume some familiarity in this paper of modern algebra. In particular, of polynomial rings and quotients. However, large swaths of this paper should be intelligible to a motivated high school student, so long as they gloss over some of the proofs. Our main audience, is thus undergraduates interested in research opportunities related to this work. Large parts of this paper, in fact, were written while the author was an undergraduate, under the supervision of Dr. James Cook at Liberty University. 

In some sense then, this paper is home to some of the miscellaneous results that we discovered during our summer research session, which were eventually found to be not directly relevant to the other papers we were writing at the time, but are nevertheless important to the general context of our research project, which involves seeing how much of the typical undergraduate calculus sequence (Calculus I, Calculus II, Differential Equations, etc...) can be developed over an algebra. Thus setting the context of why one would write such a paper in the first place, which is taking our original research program\footnote{For a general introduction to our research program, see Cook \cite{Cook12}.} backwards a step from calculus to algebra, in the hopes of gaining a broader insight into our work.  

\section{Introduction}

 Algebras are an interesting structure to study for undergraduates because they build on the already familiar structure of a real vector space:

\begin{definition}
	An algebra $\mathcal{A}$ is a finite dimensional real vector space together with a bilinear multiplication operation $\star : \mathcal{A} \times \mathcal{A} \rightarrow \mathcal{A}$ satisfying the following properties:
	\begin{enumerate}
		\item $ v \star (w \star z) = (v \star w) \star z$ for all $v,w,z \in \mathcal{A}$
		\item There exists an element $\mathbbm{1} \in \mathcal{A}$ such that $\mathbbm{1} \star z = z \star \mathbbm{1} = z$ for all $z \in \mathcal{A}$.
	\end{enumerate}
\end{definition}

 Making use of this vector space structure, many algebraic results may be derived from taking advantage of this structure, in particular, from the fact that every finite dimensional vector space has a basis. In particular, the following result is useful: 

\begin{proposition}
	\label{prop:homo_basis_elements}
	Given an algebra $\mathcal{A}$ with basis $\beta = \{ v_1, \dots , v_n \}$ and a linear map $\phi : \mathcal{A} \rightarrow \mathcal{B}$ between two algebras, if $\phi(v_i \star v_j) = \phi(v_i) \star \phi(v_j)$ for all basis elements $v_i, v_j$ and $\phi(\mathbbm{1}) = \mathbbm{1}$, then $\phi(v \star w) = \phi(v) \star \phi(w)$ for all $v,w \in \mathcal{A}$.
\end{proposition}

 Stated in more abstract terms, the fact that our algebras are \textit{finite dimensional} may be characterized by the following algebraic conditions:

\begin{proposition}
	\label{prop:algebra_noetherian}
	Let $\mathcal{A}$ be a commutative algebra, then $\mathcal{A}$ is a Noetherian and Artinian ring, In particular, this implies that every ideal $I$ of $\mathcal{A}$ is finitely generated.
\end{proposition}

 Although, as mentioned, many algebraic results follow simply from exploiting the vector space structure of an algebra, this proposition will be important when connecting our work to the more general algebraic context of ring theory.

Beyond conceptual familiarity, however, the vector space view of algebras gives us another important tool: matrix representations.

\begin{definition}
	Given an algebra $\mathcal{A}$ with basis $\beta = \{ v_1, \dots, v_n \}$, and a fixed element $\alpha \in \mathcal{A}$, notice that the map $L_\alpha(z) = \alpha \star z$ is a linear transformation from $\mathcal{A}$ to itself. We claim that under the standard addition and scalar multiplication of linear maps, treating composition of linear maps as multiplication that the collection of these linear maps forms an algebra isomorphic to $\mathcal{A}$.
	
	In particular, if we consider the standard matrices of these linear transformations with respect to the basis beta, defining $M_\beta(\alpha) = [L_\alpha]_\beta$, then these too form an algebra isomorphic to $\mathcal{A}$, which we call the regular representation of the algebra with respect to the basis $\beta$, and denote by $M_\beta(\mathcal{A})$. Furthermore, $M_\beta : \mathcal{A} \rightarrow M_\beta(\mathcal{A})$ gives an isomorphism between these algebras.\footnote{For a full presentation of this, and proof of the theorems we implicitly use in this definition, see \cite{Cook12}.}
	
\end{definition}

\noindent Below we give an example of how to calculate $M_\beta(z)$ for an arbitrary element $z \in \mathcal{A}$.

\begin{example}\label{ex:2}
	Let $\beta = \{ 1, i \}$ and let $x = x + i y$ be an arbitrary element of $\mathbb{C}$. Applying the same method as we did in the preceding example, we find $(x + i y) \star 1 = x + i y$ and $(x + i y) \star i = -y + i x$. Hence: $$ M_\beta(z) = \begin{bmatrix} x & -y \\ y & x \end{bmatrix} $$
\end{example}

 Matrix representations were used, for instance, in \cite{Cook12}, in order to define the notion of an $\Acal$-derivative, as well as in \cite{Freese15} in order to prove the k-Pythagorean theorem, a generalization of the usual identity $\sin^2(z) + cos^2(z) = 1$. 
	
\section{Algebras and Algebra Presentations}
\label{sec:Algebra_Presentations}

Another convenient way to describe a finite dimensional, commutative, and unital algebra is as the quotient of some real polynomial ring by an ideal. If we have such an algebra $\mathcal{A}$ isomorphic to $\mathcal{P} = \mathbb{R}[x_1, \dots , x_k]/I$ for some $k \in \mathbb{N}$ and $I$ an ideal of $\mathbb{R}[x_1, \dots, x_k]$ we say that $\mathcal{P}$ is a \textit{presentation} of the algebra $\mathcal{A}$. While this moves us somewhat further away from the more familiar context of finite dimensional vector spaces, it gives us an intuitive view of algebras in terms of generators and relations which will connect our program with our original motivations (generalizing the construction by which $\mathbb{C}$ is constructed from $\mathbb{R}$ by ``adjoining $i$'' and imposing the relation $i^2 = -1$).

\begin{definition}[Standard presentations of typical algebras] $ $
	\label{def:common_algebras}
	\begin{enumerate}
		\item The $n$-hyperbolic numbers: $\mathcal{H}_n := \mathbb{R}[j]/\langle j^n - 1 \rangle$
		\item The $n$-complicated numbers: $\mathcal{C}_n := \mathbb{R}[i]/\langle i^n + 1 \rangle$
		\item The $n$-nil numbers: $ \mathbf{\Gamma}_n := \mathbb{R}[\epsilon]/\langle \epsilon^n \rangle$
		\item The total $n$-nil numbers: $\mathbf{\Xi}_n := \mathbb{R}[\epsilon_1,\dots,\epsilon_n]/\langle \epsilon_i\epsilon_j | i,j \in \{ 1, 2, \dots , n \} \rangle$
	\end{enumerate}
\end{definition}

For example, $\mathcal{C}_2$ is just the usual complex numbers, denoted simply $\mathbb{C}$, and $\mathcal{H}_2$ is just the hyperbolic numbers, denoted simply $\mathcal{H}$. Similarly, we take the convention that $\mathbf{\Gamma}$ by itself denotes $\mathbf{\Gamma}_2$. 

The nil numbers are a special class of what in our terminology we will call \textit{unital nil algebras} -- that is, an algebra with basis $\{1, \epsilon_1, \dots, \epsilon_{n-1} \}$, where each $\epsilon_k$ is \textit{nilpotent}. In other words, for each $\epsilon_k$ there exists $m \in \mathbb{N}$ such that $(\epsilon_k)^m = 0$. This terminology is inspired by the use of  the term \textit{nil algebra} used by Abian \cite{abian71} to refer to an algebra in which every element of the algebra is nilpotent. Although, we will also sometimes use ``nil algebra'' instead of ``unital nil algebra'' here for brevity, as in our context all algebras are unital. A unital nil algebra is in some sense the closest you can get to a nil algebra while still being a unital algebra. 

\begin{definition}
	\label{def:nil_star}
	Given an algebra $\mathcal{A}$, let $\mathrm{Nil}^*(\mathcal{A})$ denote the smallest unital sub-algebra of $\mathcal{A}$ that contains the nilradical of $\mathcal{A}$, $\mathrm{Nil}(\mathcal{A})$ -- that is, the ideal formed from all nilpotent elements of $\mathcal{A}$.
\end{definition}

\begin{proposition}
	An algebra $\mathcal{A}$ is a unital nil algebra if and only if it is the smallest unital subalgebra of $\mathcal{A}$ containing the nilradical of $\mathcal{A}$. In other words, $\mathcal{A}$ is a unital nil algebra if and only if $\mathcal{A} = \mathrm{Nil}^*(\mathcal{A})$.
	
	\begin{proof}
		If $\mathcal{A}$ is a unital nil algebra with unital nil basis $\{1, \epsilon_1, \dots, \epsilon_n \}$, then clearly $\mathrm{Nil}(\mathcal{A}) = \langle \epsilon_1, \dots, \epsilon_n \rangle$, and $\mathcal{A}/\mathrm{Nil}(\mathcal{A}) \cong \mathbb{R}$, so $\mathrm{Nil}(\mathcal{A})$ is maximal, and hence $\mathcal{A}$ is the smallest unital nil algebra containing $\mathrm{Nil}(\mathcal{A})$.
		
		Conversely, suppose that $\mathcal{A} = \mathrm{Nil}^*(\mathcal{A})$, and let $\{ \epsilon_1, \epsilon_2, \dots, \epsilon_n \}$ be a basis for $\mathrm{Nil}(\mathcal{A})$, then clearly $\{ 1, \epsilon_1, \dots, \epsilon_n \}$ is a linearly independent set that spans a subset of $\mathrm{Nil}^*(\mathcal{A})$. Also, we argue that this set must span $\mathrm{Nil}^*(\mathcal{A})$, since it contains all of $\mathrm{Nil}(\mathcal{A})$, and the unit $1$, so by the minimality condition, this must be all of $\mathrm{Nil}^*(\mathcal{A})$. Thus, $\mathcal{A} = \mathrm{Nil}^*(\mathcal{A})$ is a multiplicative nil algebra.
	\end{proof}
\end{proposition}

In addition to these basic families, we should also mention the so called $n$-complex numbers $\mathbb{C}_n = \mathbb{C}^{\otimes n}$, where $X^{\otimes n}$ denotes the $n$-fold tensor product of rings.\footnote{For the unfamiliar reader, the tensor product of algebras can be thought of as the algebra which combines the set of generators and relations for a presentation for the algebra. In other words, $\mathbb{R}[x_1, \dots, x_n]/I \otimes \mathbb{R}[y_1, \dots, y_m]/J \cong \mathbb{R}[x_1,\dots,x_n,y_1,\dots,y_m]/(I+J)$. So $\mathbb{C}_2 = \mathbb{C} \otimes \mathbb{C} \cong \mathbb{R}[i_1,i_2]/\langle i_1^2 + 1, i_2^2 +1 \rangle$ for example.} In particular, the analysis of the bicomplex numbers $\mathbb{C}_2$ has been studied extensively, for example, by Price \cite{Price}.

Certain presentations of an algebra are more economical than others. For example, in Definition \ref{def:common_algebras} we defined the algebra $\mathcal{H}$ by the presentation $\mathbb{R}[j]/\langle j^2 -1 \rangle$, but $\mathcal{H}$ could also be presented as $\mathbb{R}[j,k]/\langle j^2 - 1, k \rangle$. This motivates the following definition:

\begin{definition}
	Given a presentation $\mathcal{P} = \mathbb{R}[x_1,\dots,x_k]/I$ of an algebra $\mathcal{A}$, we say that the presentation $\mathcal{P}$ is \textit{degenerate} if the set $\{ 1 + I, x_1 + I, \dots , x_k + I \}$ is linearly dependent as a vector space over $\mathbb{R}$. 
\end{definition}

From this, it is easy to see that the standard presentation for $\mathcal{H}$ is non-degenerate, but that $\mathbb{R}[j,k]/\langle j^2 - 1, k \rangle$ is degenerate.

It also also oftentimes convenient to identify the elements $1,x_1, \dots, x_n$ of a presentation with a basis for the algebra being presented, motivating yet another definition:

\begin{definition}
	Let $\mathcal{P} = \mathbb{R}[x_1,\dots,x_n]/I$ be a presentation of the algebra $\mathcal{A}$, then we say $\mathcal{P}$ is a \textit{basic presentation} if the set $\{ 1 + I, x_1 + I, \dots , x_n + I \}$ forms a basis for $\mathcal{P}$ with respect to the real vector space structure on $\Pcal$.
\end{definition}

Most of the time a basic presentation is not the most economical way to describe an algebra. For example, the standard presentation of $\mathcal{H}_3$, $\mathbb{R}[j]/\langle j^3 - 1 \rangle$ is not a basic presentation because $\{ 1, j \}$ does not form a basis for $\mathcal{H}_3$. However, this presentation is simpler than the basic presentation $\mathbb{R}[x,y]/\langle x^2 - y, y^2 - x, xy - 1 \rangle$ of $\mathcal{H}_3$. One of the notable exceptions to this rule might be the class of totally nil numbers $\mathbf{\Xi}_n$, which the reader should confirm was defined in Definition \ref{def:common_algebras} using only basic presentations. However, this kind of presentation has the advantage of connecting the vector space viewpoint of algebras with the generators and relations viewpoint\footnote{In an unpublished draft of our paper ``Logarithms Over a Real Associative Algebra'', we used this notion of a basic presentation to give an elementary proof of a classification theorem for commutative algebras. However, for brevity's sake, this was later revised to take advantage of the standard classification of Artinian rings into a direct product of local rings.}.

To fully understand this connection between the vector space and the generators and relations conceptions of algebras, we need to introduce yet another notion. If $\Acal$ is a real vector space with basis $\beta = \{ v_1, \dots , v_n \}$ then given appropriate \textit{structure constants} $c_{ij}^k \in \mathbb{R}$ we may define a multiplication on $\Acal$. In particular, define
$$ v_i \star v_j = \sum_{k=1}^n c_{ij}^k v_k $$ on basis elements, and extend bilinearly to define $\star$ on $\Acal$. Naturally, the structure constants must be given such that the defined multiplication is associative and unital. That said, we typically begin with a given algebra $\Acal$ and simply use the structure constants with respect to a given basis to study the structure of $\Acal$. For example:

\begin{theorem}
	\label{thm:exists_basic_presentation}
	Given any algebra $\mathcal{A}$ with structure constants $c_{ij}^k$ associated with the basis $\beta = \{v_1, \dots, v_n \}$ where $v_1=1$, the presentation $\mathcal{P}$ defined by: $$ \mathbb{R}[\bar v_1, \dots, \bar v_n]/\langle \{\bar v_i \bar v_j - \sum_{k=1}^n c_{ij}^k \bar v_k | i, j = 2, \dots, n \} \rangle$$ is a basic presentation of $\mathcal{A}$.
	
	\begin{proof}
		First, let us prove that $\mathcal{P}$ is indeed a presentation of $\mathcal{A}$. We must show that $\phi : \mathcal{A} \rightarrow \mathcal{P}$ defined by setting $\phi(v_i) = \bar v_i + I$ on basis elements on extending linearly forms a linear bijection. Since we defined $\phi$ by linear extension, we already know that $\phi$ is linear, so it remains to show bijectivity and the homomorphism property. 
		
		By proposition \ref{prop:homo_basis_elements}, it suffices to show the homomorphism property on basis elements. Thus, making use of the fact that the presentation of $\mathcal{P}$ forces relations of the form $$ \bar v_i \bar v_j = \sum_{k=1}^n c_{ij}^k \bar v_k $$ to hold in the algebra, we consider:
		\begin{equation*}
		\begin{aligned}
		\phi(v_i \star v_j) &= \phi \left(\sum_{k=1}^n c_{ij}^k v_k \right) = \sum_{k=1}^n c_{ij}^k \phi(v_k) \\ &= \sum_{k=1}^n c_{ij}^k (\bar v_k+I) = \left(\sum_{k=1}^n c_{ij}^k \bar v_k \right) + I \\ &= \bar v_i \bar v_j + I = \phi(v_i)\phi(v_j)
		\end{aligned}
		\end{equation*}
		
		 \noindent And hence, showing that $\phi$ is an algebra homomorphism.
		
		To show that $\phi$ is injective, notice that $\phi(c_1v_1 + \dots c_nv_n) = \phi(d_1v_1 + \dots + d_nv_n)$ if and only if $(c_1 - d_1) \bar v_1 + \dots + (c_n - d_n) \bar v_n = I$, and hence if and only if $c_i = d_i$ for all $i$, thus showing $c_1 v_1 + \dots + c_n v_n = d_1 v_1 + \dots + d_n v_n$. (Since the only non-zero elements of $I$ are second order) Also, $\phi$ is surjective since we can use the relation $$ \bar v_i \bar v_j = \sum_{k=1}^n c_{ij}^k \bar v_k $$ induced by the quotient to successively reduce polynomials in $\bar v_1, \bar v_2, \dots , \bar v_n$ into elements of the form $\phi(c_1 v_1 + \dots + c_n v_n) = c_1 \bar v_1 + \dots + c_n \bar v_n$.
		
		Finally, since $\phi$ is a linear bijection, by linear algebra we know that since $\beta$ forms a basis for $\mathcal{A}$, $\phi(\beta) = \{ \phi(v_1), \dots , \phi(v_n) \} = \{ \bar v_1 + I, \bar v_2 + I, \dots ,\bar v_n + I \}$ is a basis for $\phi(\mathcal{A}) = \mathcal{P}$, and hence by definition $\mathcal{P}$ is a basic presentation of $\mathcal{A}$.
	\end{proof}
\end{theorem}

We call the presentation constructed in the previous theorem the \textit{canonical basic presentation} with respect to the basis $\beta = \{ v_1, \dots , v_n \}$. Combining Theorem \ref{thm:exists_basic_presentation} and the fact that each basic presentation is non-degenerate we obtain the following reassuring corollary:

\begin{corollary}
	Every algebra $\mathcal{A}$ has a non-degenerate presentation $\mathcal{P}$.
\end{corollary}

It is also useful to note when proving propositions involving algebra presentations that instead of explicitly showing the homomorphism property for a linear bijection to show that it is an algebra isomorphism, we can make an argument more directly in terms of the generators and relations defining the presentation, but to formalize this we must first make some new definitions:

\begin{definition}
	Given an algebra $\mathcal{A}$ with canonical basic presentation $\mathbb{R}[v_1, \dots, v_n]/I$ and isomorphism $\phi : \mathcal{A} \rightarrow \mathbb{R}[\bar v_1, \dots, \bar v_n]/I$ as given in the proof of Theorem \ref{thm:exists_basic_presentation}, let $\psi_I : \mathbb{R}[v_1, \dots, v_n] \rightarrow \mathbb{R}[v_1, \dots, v_n]/I$ be the natural quotient map, then we define the evaluation homomorphism of the algebra with respect to the basic presentation $\mathrm{ev}_\mathcal{A} : \mathbb{R}[v_1, \dots, v_n] \rightarrow \mathcal{A}$ by $\mathrm{ev}_\mathcal{A} = \phi^{-1} \circ \psi_I $.
\end{definition}

The idea is that given a polynomial in the basis elements of the algebra -- representing a formal expression in the algebra, we can evaluate that expression to produce an algebra element. For example, consider the basic presentation $\mathbb{R}[\hat j]/\langle \hat j^2 -1 \rangle$ for $\mathcal{H}$ with the natural isomorphism between the two algebras sending $\hat j$ to $j$, then $\hat j^2 + \hat j - 1 \neq \hat j$, since these are just formal polynomials in $\mathbb{R}[\hat j]$, but $\mathrm{ev}_\mathcal{H}(\hat j^2 + \hat j - 1) = j^2 + j - 1 = j = \mathrm{ev}_\mathcal{H}(\hat j)$. From this, we can represent the notion of two algebras satisfying the same relations precisely using a commutative diagram, as we do in the following theorem:

\begin{theorem}
	\label{thm:generator_relns}
	Given a bijective linear map $\psi : \mathcal{A} \rightarrow \mathcal{B}$ between two algebras and a basis $\beta = \{ v_1, \dots, v_n \}$ of $\mathcal{A}$, consider the canonical basic presentations of $\mathcal{A}$ with respect to the basis $\beta$, and $\mathcal{B}$ with respect to the basis $\psi(\beta)$, then let $\overline\psi : \mathbb{R}[\overline v_1,\dots,\overline v_n] \rightarrow \mathbb{R}[\overline{\psi(v_1)}, \dots,\overline{\psi(v_n)}]$ be the natural isomorphism defined by setting $\overline\psi(v_i) = \overline{\psi(v_i)}$ and extending linearly. If the following diagram commutes, then $\psi$ is an algebra isomorphism: $$
	\begin{CD}
	\mathbb{R}[\overline v_1, \dots \overline v_n]    @>\overline{\psi}>> \mathbb{R}[\overline{\psi(v_1)}, \dots, \overline{\psi(v_n)}] \\
	@VV \mathrm{ev}_\mathcal{A} V        @VV \mathrm{ev}_\mathcal{B} V \\
	\mathcal{A}                 @>  \psi  >>      \mathcal{B}
	\end{CD}
	$$
	
	\begin{proof}
		Since we already know that $\psi$ is a linear bijection, it suffices to prove the homomorphism property on basis elements $v_1, \dots, v_n$. Notice that the given diagram commuting is equivalent to the statement that for all formal polynomials $p(v_1,\dots,v_n) \in \mathbb{R}[v_, \dots, v_n]$ $$ \psi(\mathrm{ev}_\mathcal{A}(p(\overline v_1, \dots, \overline v_n))) = \mathrm{ev}_\mathcal{B}(p(\overline{\psi( v_1)},\dots,\overline{\psi(v_n)}): $$ Or, since $\mathrm{ev}_\mathcal{A}$ and $\mathrm{ev}_\mathcal{B}$ are algebra homomorphisms: $$ \psi(p(\mathrm{ev}_\mathcal{A}(\overline v_1),\dots,\mathrm{ev}_\mathcal{A}(\overline v_n))) = p(\mathrm{ev}_\mathcal{B}(\overline{\psi(v_1)}),\dots,\mathrm{ev}_\mathcal{B}(\overline{\psi(v_n)})) $$ Hence, applying this equality to $p = v_i \star v_j$ implies $$ \psi(v_i \star v_j) = \psi(v_i) \star \psi(v_j) $$ since under the canonical basic presentation of $\mathcal{A}$ we have $\mathrm{ev}_\mathcal{A}(\overline v_i) = v_i$ and $\mathrm{ev}_\mathcal{B}(\overline{\psi(v_j)}) = \psi(v_j)$. Hence, we have shown that the homomorphism property holds for all basis elements, and therefore that $\psi$ is an algebra isomorphism.
	\end{proof}
\end{theorem}

 It is important to note however that intuitively what Theorem \ref{thm:generator_relns} says is that the algebra $\mathcal{A}$ satisfies a relation $p(v_1, \dots, v_n) = 0$ between its basis elements $v_1, \dots, v_n$  if and only if the corresponding relation $p(\psi(v_1),\dots,\psi(v_n)) = 0$ holds in $\mathcal{B}$. The commutative diagram is simply a concise way to formalize this idea.
	
	\section{Polynomials and Irreducibility Over an Algebra}
	
	Unlike in the case of a field, the notion of degree for polynomials in an algebra requires a certain amount of care. For example, in $\mathcal{H}$, $2jz + 1$ factors as $((j-1)z + 1)((j+1)z+1)$. Thus, in general we will not have $deg(f(z)g(z)) = deg(f(z)) + deg(g(z))$. However, we still may still use the standard definition for polynomial rings $R[x]$.
	
	\begin{definition}[]
		let $f(z) \in \mathcal{A}[z]$, then by definition of a polynomial ring, $f(z) = a_0 + a_1 z + a_2 z^2 + \dots $ for some $a_0,a_1, \dots \in \mathcal{A}$ where only finitely many $a_i$ are non-zero. We define the degree of $f(z)$ to be $\deg(f(z)) = n$, where n is the largest integer such that $a_n \neq 0$. 
	\end{definition}

	 \noindent We also still have a notion of irreducibility for polynomials over an algebra:
	
	\begin{definition}[Irreducibility]
		Let $f(z) \in \mathcal{A}[z]$, we say that $f(z)$ is \textit{irreducible} over $\mathcal{A}$ if $f(z) = g(z)h(z)$ for some $g(z),h(z) \in \mathcal{A}[z]$ implies either $g(z)$, or $h(z)$ must be of degree zero.
	\end{definition}
	
	\begin{example}
		\label{ex:irred_poly}
		Consider the polynomial $z^2 + jz + j \in \mathcal{H}[z]$. If $z^2 + j z + j$ is reducible, then since it is a second order monic polynomial, $z^2 + j z + j = (z + a + b j)(z + c + dj)$ for some $a,b,c,d \in \mathbb{R}$. Equating coefficients then, it can be shown by explicit calculations that the resulting system of equations has no real solutions, and consequently that the polynomial $z^2 + j z + j$ is irreducible.
	\end{example}

	 This notion is important for the study of $\Acal$-ODEs, where it turns out that constant coefficient differential operators over an algebra $\Acal[D]$ have the same structure as polynomials over the algebra.
	
	Although usually proven for an integral domain, it is also important to note that the division algorithm still holds for the algebras that we consider in this paper, as is shown in McCoy's book \cite{mccoy}:
	
	\begin{theorem}
		Let $R$ be an arbitrary ring with unit element, and let $f(x) = a_n x^n + \dots + a_0, \; g(x) = b_m x^m + \dots + b_0$ where $a_n, b_m \neq 0$ and $b_m$ is a unit in $R$, then there exist unique elements $q(x),p(x),r(x),s(x) \in R[x]$ such that: $$ f(x) = q(x)g(x) + r(x) \ \ \& \ \ f(x) = g(x)p(x) + s(x)$$ where $r(x)$ and $s(x)$ are either $0$ or of degree less than $m$, and $q(x),p(x)$ are either both zero or of degree $n-m \geq 0$
	\end{theorem}
	
	 \noindent As a corollary of this, there is also a factor theorem in our context:
	
	\begin{theorem}
		\label{thm:factor_theorem_diff}
		Let $\mathcal{A}$ be an associative commutative algebra, and let $f(z) \in \mathcal{A}[z]$, then $\alpha \in \mathcal{A}$. $f(\alpha) = 0$ if and only if $f(x) = (x-\alpha)g(x)$ for some $g(x) \in \mathcal{A}[z]$
	\end{theorem}
	
	\section{Semisimple Algebras}
	
	An important notion throughout this paper is that of a semisimple algebra. Usually, an algebra $\mathcal{A}$ is  defined to be semisimple if and only if its Jacobson radical is trivial, but fortunately in our context, the Jacobson radical coincides with the conceptually simpler \textit{nilradical}, which recalling Definition \ref{def:nil_star} is the set of all elements $z \in \mathcal{A}$ such that $z^n = 0$ for some $n \in \mathbb{N}$, which we denote by $\mathrm{Nil}(\mathcal{A})$
	
	\begin{theorem}
		In a finite dimensional associative algebra $\mathcal{A}$, the Jacobson radical and nilradical coincide, and hence, a finite dimensional associative algebra $\mathcal{A}$ is semisimple if and only if its only nilpotent element is $0$.
		
		\begin{proof}
			Since $\mathcal{A}$ is a finite dimensional algebra, by Proposition \ref{prop:algebra_noetherian} it is Artinian, and hence the maximal ideals of $\mathcal{A}$ coincide with the prime ideals of $\mathcal{A}$. Therefore, since the $Nil(\mathcal{A})$ is the intersection of all the prime ideals of $\mathcal{A}$, and the Jacobson radical $J(\mathcal{A})$ is the intersection of all maximal right ideals of $\mathcal{A}$, $J(\mathcal{A}) = \mathrm{Nil}(\mathcal{A})$.
		\end{proof}
	\end{theorem}
	
	 One of the main reasons why it is often nicer to work in a semisimple algebra is because of the well-known classification of such rings given by Artin and Wedderburn:
	
	\begin{theorem}[Artin-Wedderburn Theorem]
		\label{Artin-Wedderburn}
		If $R$ is an an Artinian semisimple ring, then $R$ is isomorphic to a product of finitely many matrix rings over division algebras.
	\end{theorem}
	
	 For a more in depth discussion and proof of the Artin-Wedderburn theorem, the reader may consult \cite{dummit03}, or \cite{abian71} for a discussion of the theorem in the context of associative algebras. In particular, for us this theorem means that any finite dimensional real associative algebra $\mathcal{A}$ will be isomorphic to a finite product of matrix rings over $\mathbb{R}$, $\mathbb{C}$, or $\mathbb{H}$, as these are the only finite dimensional associative division algebras over the reals. In addition to this, if we restrict ourselves to the case of a commutative algebra, we obtain:
	
	\begin{corollary}
		\label{cor:commutative_wedderburn}
		Every finite dimensional semisimple commutative algebra is isomorphic to the direct product of $m$ copies of $\mathbb{R}$ and $k$ copies of $\mathbb{C}$.
	\end{corollary}

	\noindent As $\mathbb{H}$, and any $n$ by $n$ matrix ring over $\mathbb{R}, \mathbb{C}$, or $\mathbb{H}$ will be non-commutative for $n > 1$.
	
	\begin{example}
		\label{ex:H_iso}
		$\mathcal{H} \cong \mathbb{R} \times \mathbb{R}$, which may be seen from the explicit isomorphism $\phi : \mathbb{R} \times \mathbb{R} \rightarrow \mathcal{H}$ defined by $\phi(x,y) = \frac{1}{2}(x+y) + \frac{j}{2}(x-y)$
	\end{example}
	
	Clearly, $\mathbf{\Gamma}_n$ and $\mathbf{\Xi}_n$ are not semisimple algebras, but the other families of algebras we defined in Section \ref{sec:Algebra_Presentations}, namely $\mathcal{H}_n, \mathcal{C}_n,$ and $\mathbb{C}_n$ are, and so by the Artin-Wedderburn Theorem are all isomorphic to a direct product of $m$ copies of $\mathbb{R}$ and $k$ copies of $\mathbb{C}$. To better understand their structure, we wish then to determine precisely how many copies of $\mathbb{R}$ and $\mathbb{C}$ each of these algebras contain in their Wedderburn decomposition. To accomplish this, we must first recall a theorem from algebra:
	
	\begin{theorem}[\cite{dummit03}]
		\label{thm:chinese_remainder}
		Let $R$ be a ring and $A_1,\dots,A_n$ be ideals of $R$, then if the ideals $A_i,A_j$ are comaximal for all $i,j = 1, \dots, n$ we have the isomorphism: $$R/(A_1A_2 \dots A_n) \cong R/A_1 \times R/A_2 \times \dots \times R/A_n$$
	\end{theorem}
	
	 Finally, we also need the following proposition, which uses the structure of $n$th roots of $1$ and $-1$ in $\mathbb{C}$ to deduce the factorizations of $x^n - 1$ and $x^n + 1$ into irreducible polynomials over $\mathbb{R}$.
	
	\begin{proposition}
		\label{prop:roots_of_unity}
		Over $\mathbb{C}$ consider the equation $z^n + 1 = 0$, these are the $n$th roots of $-1$. If $n$ is even, then all $n$ roots of this equation are complex, and hence the polynomial $x^n - 1$ may be factored as the product of $n/2$ irreducible quadratic polynomials over the reals. If $n$ is odd, then we know $z^n + 1 = 0$ has one real root and $n-1$ complex roots coming in conjugate pairs over $\mathcal{C}$. Hence, $x^n - 1$ factors as the product of one real linear factor and $(n-1)/2$ irreducible quadratics over $\mathbb{R}$.
		
		Similarly, if we consider $z^n - 1 = 0$, if $n$ is even then the equation has two real roots $z = 1,-1$, and $n-2$ complex roots coming in conjugate pairs. Hence, $x^n - 1$ factors over the reals as the product of two real linear factors, and $(n-2)/2$ irreducible quadratic factors. If $n$ is odd, then $z^n - 1 = 0$ has only a singe real root $z = 1$, and $n-2$ remaining complex roots coming in conjugate pairs. Thus, $x^n - 1$ factors as the product of a single real linear factor, and $(n-1)/2$ irreducible quadratics. 
	\end{proposition}
	
	 \noindent Hence, we obtain the following:
	
	\begin{corollary}
		For all $k \in \mathbb{Z}^+$ we have the following isomorphisms:
		\begin{enumerate}
			\item $\mathcal{H}_{2k} \cong \mathbb{R}^2 \times \mathbb{C}^{k-1}$
			\item $\mathcal{C}_{2k} \cong \mathbb{C}_k \cong \mathbb{C}^k$
			\item $\mathcal{H}_{2k + 1} \cong \mathcal{C}_{2k + 1} \cong \mathbb{R} \times \mathbb{C}^k$
		\end{enumerate}
		
		\begin{proof}
			Consider $\mathcal{H}_{2k} = \mathbb{R}[j]/\langle j^{2k} - 1\rangle$. By Proposition \ref{prop:roots_of_unity}, we have $j^{2k} - 1 = (j + a_1)(j + a_2)(j^2 + b_1 j + c_1)\dots(j^2 + b_{k-1}j + c_{k-1})$ as a factorization of $j^{2k}-1$ into irreducible polynomials. Also, clearly each of these factors are coprime, and hence since $\mathbb{R}[j]$ is a PID, the ideals generated by each of the factors is comaximal, and thus Theorem \ref{thm:chinese_remainder} applies, so we obtain:
			\begin{equation*}
			\begin{aligned}
			\mathcal{H}_{2k} &= \mathbb{R}[j]/\langle j^{2k} - 1\rangle \\ &= \mathbb{R}[j]/\langle (j + a_1)(j + a_2)(j^2 + b_1 j + c_1)\dots(j^2 + b_{k-1} j + c_{k-1}) \rangle \\ &\cong \mathbb{R}[j]/\langle j + a_1 \rangle \langle j + a_2 \rangle \langle j^2 + b_1 j + c_1 \rangle \dots \langle j^2 + b_{k-1}j + c_{k-1} \rangle  \\ &\cong \mathbb{R}[j]/\langle j + a_1 \rangle \times \mathbb{R}[j]/\langle j + a_2 \rangle \times \mathbb{R}[j]/\langle j^2 + b_1 j + c_1 \rangle \times \dots \\ 
			& \qquad \dots \times \mathbb{R}[j]/\langle j^2 + b_{k-1} j + c_{k-1} \rangle \\ &\cong \mathbb{R} \times \mathbb{R} \times \mathbb{C} \times \dots \times \mathbb{C} \\ &\cong \mathbb{R}^2 \times \mathbb{C}^{k-1}
			\end{aligned}
			\end{equation*}
			
			 $(2)$ and $(3)$ may be derived similarly, applying the relevant remarks made in Proposition \ref{prop:roots_of_unity} with the exception of the isomorphism $\mathbb{C}_{2k} \cong \mathbb{C}^k$.
			
			To show $\mathbb{C}_{k} \cong \mathbb{C}^k$, note that since $\mathbb{C}_{k}$ is semisimple $\mathbb{C}_{k} \cong \mathbb{R}^n \times \mathbb{C}^m$ for some $n,m \in \mathbb{N}$ by the Wedderburn decomposition. Also, note that $\mathbb{C}_{k}$ contains an element $\mathbf{i}$ such that $\mathbf{i}^2 = -1$, so the same must hold in the Wedderburn decomposition of $\mathbb{C}_{k}$. Hence, we must have some element $(x_1,\dots,x_n,z_1,\dots,z_m) \in \mathbb{R}^n \times \mathbb{C}^m$ such that $(x_1,\dots,x_n,z_1,\dots,z_m)^2 = (x_1^2,\dots,x_n^2,z_1^2,\dots,z_m^2) = (-1,-1,\dots,-1)$, which clearly is impossible unless $n = 0$, and hence, the only option for the Wedderburn decomposition is $\mathbb{C}^k$. Therefore, by elimination it must be the case that $\mathbb{C}_k \cong \mathbb{C}^k$.
			
		\end{proof}
	\end{corollary}
	
\section{Structure of Zero Divisors in an Algebra}

It was shown by Freese \cite{Freese15} that zero divisors play a significant role in the study of analysis over an algebra. In future work we will show that the presence of zero divisors is similarly pervasive in our overall study. Thus, we wish to characterize the zero divisors of algebras. Moreover, to better understand the structure of how said zero divisors function in an commutative algebra, we would also like to be able to describe the \textit{annihilators} of the zero divisors in an algebra\footnote{Notice that for the non-commutative case, we must consider both right and left annihilators. Our future work on logarithms hold only in the commutative case, and thus we do not consider the non-commutative case here. However, it is likely that some of our remarks could be generalized.} That is, given a zero divisor $a \in \mathcal{A}$ we would like to understand the set of all elements $x \in \mathcal{A}$ such that $xa = 0$, which we denote by $\mathrm{Ann}(a)$. 

\begin{definition}
	Given a commutative algebra $\Acal$, we denote the set of zero divisors in the algebra $\zd(\Acal) = \{ z \in \Acal \; | \; z \star w = 0, \text{for some } w \in \Acal, w \neq 0 \}$. We also let $\zd^* = \zd - \{ 0 \}$ denote the set of non-trivial zero divisors in $\Acal$.
\end{definition}

\noindent To begin our exploration into the structure of the zero divisors we note that: 

\begin{theorem}
	\label{thm:zero_divisor_regular_rep}
	An element $a \in \mathcal{A}$ is a zero divisor if and only if $\det{M_\beta(a)} = 0$.
\end{theorem}

 One consequence of this theorem is the fact that units are dense in an algebra, since the equation specifying zero divisors in an $n$-dimensional algebra will have at most a solution set which is $n-1$ dimensional.

Corollary \ref{cor:commutative_wedderburn}, gives us an important result that lets us characterize the zero divisors of a semisimple algebra:

\begin{theorem}
	\label{thm:zero_divisors_semisimple}
	Let $\mathcal{A}$ be a semisimple commutative algebra, with $\mathcal{A} \cong \mathbb{R}^{m} \times \mathbb{C}^{k}$ for some $m,k$, so let $\phi : \mathcal{A} \rightarrow \mathbb{R}^{m} \times \mathbb{C}^{k}$ be an isomorphism. Then $a \in \zd(\Acal)$ if and only if $\phi(a)$ is zero in at least one component. Moreover, the annihilators of $a$ are exactly the elements $b \in \mathcal{A}$ corresponding to $\phi(b) \in \mathbb{R}^{m} \times \mathbb{C}^{k}$ with zeros in the components which were non-zero in $\phi(a)$.
	
	\begin{proof}
		Algebra isomorphisms are ring isomorphisms, and hence preserve zero divisors and any other relevant properties of zero divisors, such as their annihilators. Hence, since $\mathcal{A} \cong \mathbb{R}^{m} \times \mathbb{C}^{k}$, we simply need to characterize the zero divisors and their annihilators in $\mathbb{R}^{m} \times \mathbb{C}^{k}$.
		
		\noindent Let $a = (a_1, \dots , a_n, b_1, \dots, b_m) \in \mathbb{R}^n \times \mathbb{C}^m$ be a zero divisor, then there exists $$b = (x_1, \dots, x_n, z_1, \dots, z_m) \in \mathbb{R}^n \times \mathbb{C}^m$$ such that 
		\begin{align*}
		(a_1, \dots , a_n, b_1, \dots, b_m)(x_1, \dots, x_n, z_1, \dots, z_m) 
		\\ = (a_1 x_1 , \dots , a_n x_n, b_1 z_1, \dots , b_m z_m) = 0 
		\end{align*}
		and thus, we must have $a_1 x_1 = \dots = a_n x_n = b_1 z_1 = \dots = b_m z_m = 0$, which implies that for all $k \in \{1, \dots , n \}$, either $a_k = 0$, or $x_k = 0$, and for all $i \in \{ 1, \dots, m \}$ either $b_i = 0$ or $ z_i = 0$. Hence, $a$ is a zero divisor, if and only if at least one of its components must be zero. Moreover, an annihilator $b$ of $a$ must have zeros in all of the components in which $a$ has no zeros.
	\end{proof}
	
\end{theorem}

\begin{example}
	Recall the isomorphism $\phi : \mathbb{R} \times \mathbb{R} \rightarrow \mathcal{H}$ defined in Example \ref{ex:H_iso}. $\phi(e_1) = \frac{1}{2} + \frac{1}{2}j$, and $\phi(e_2) = \frac{1}{2} - \frac{1}{2} j$, so by Theorem \ref{thm:zero_divisors_semisimple} we have $ \zd(\Hcal) = \mathrm{span}_\mathbb{R}\{ \frac{1}{2} + \frac{1}{2}j \} \cup \mathrm{span}_\mathbb{R} \{ \frac{1}{2} - \frac{1}{2} j \}$.
\end{example}

 The previous theorem tells us that in the commutative semisimple case the zero divisors of an algebra have a particularly nice form. We can also infer a simple geometric description of the zero divisors of a semisimple algebra from the structure of $\mathbb{R}^m \times \mathbb{C}^k$ as follows:

Consider $w = (x_1,x_2,\dots,x_m,z_1,\dots,z_k) \in \mathbb{R}^m \times \mathbb{C}^k$. By Theorem \ref{thm:zero_divisor_regular_rep} $w$ is a zero divisor if and only if $\det(M(w)) = x_1 x_2 \dots x_m |z_1|^2 |z_2|^2 \dots |z_k|^2 = 0$, which implies that at least one of $x_1,x_2, \dots, x_m$ or $z_1,z_2, \dots, z_k$ must be zero. If $x_i = 0$ for some $i$, notice there are $2k(m-1)$ remaining real degrees of freedom, and similarly if $z_i = 0$ for some $i$ there are $2m(k-1)$ remaining real degrees of freedom. Hence, the zero divisors in $\mathbb{R}^m \times \mathbb{C}^k$ consist of the union of $m$ distinct $2k(m-1)$-dimensional subspaces corresponding to the real components, and $k$ distinct $2m(k-1)$-dimensional subspaces corresponding to the complex components.

Furthermore, the annihilator of an element $z$ is simply the subspace spanned by each of the components $e_i$ where the $i$th component of $z$ is $0$. This is true since, $e_i \star e_j = 0$ for $i \neq j$.

\begin{example}
	The algebra $\mathcal{H}_3$ is generated by $j$ such that $j^3=1$. The element $z = a+jb+cj^2$ has matrix representation
	$ M(z) = \left[ \begin{array}{ccc} 
	a & c & b \\
	b & a & c \\
	c & b & a  \end{array}\right] $ for which
	$$ \text{det}(M(z)) = (a+b+c)(a^2+b^2+c^2-ab-ac-bc) $$
	It follows zero-divisors require either $a+b+c=0$ or $a^2+b^2+c^2-ab-ac-bc=0$. On the other hand, $\mathbb{R} \times \mathbb{C}$ has zero divisors $(x,0)$ and $(0,z)$ for $x \in \mathbb{R}$ and $z \in \mathbb{C}$. The isomorphism $\phi: \mathcal{H}_3 \rightarrow \mathbb{R} \times \mathbb{C}$ defined by 
	$$ \phi( a+bj+cj^2) = (a+b+c, a+be^{2\pi i/3}+ce^{4\pi i/3}). $$
	has $\phi(j) = (1,e^{2\pi i/3})$ and 
	$$ \phi^{-1}(u,x+iy) = \frac{u+2x}{3} + \frac{j}{3}\left(u-x+y\sqrt{3}\right)+\frac{j^2}{3}\left(u-x-y\sqrt{3}\right)$$ 
	Notice $\phi^{-1}(u,0) = \frac{u}{3}(1+j+j^2)$ hence identify $a=b=c=u/3$ which solves $a^2+b^2+c^2-ab-ac-bc=0$. Likewise,
	$$\phi^{-1}(0,x+iy) = \frac{2x}{3}+\frac{j}{3}\left(-x+y\sqrt{3}\right)+\frac{j^2}{3}\left(-x-y\sqrt{3}\right) =a+bj+cj^2$$ 
	provides $a=2x/3$ and $b = (-x+y\sqrt{3})/3$ and $c = (-x-y\sqrt{3})/3$ for which $a+b+c=0$.  We have shown how the zero divisors of $\mathbb{R} \times \mathbb{C}$ reveal the hidden zero divisors of $\mathcal{H}_3$.
\end{example}

Given an algebra $\mathcal{A}$ whose zero divisors we understand, a natural question to ask is:  {\it How might we characterize the zero divisors in $\mathcal{A}[z]$?} One of the most basic results in ring theory towards this end is McCoy's Theorem \cite{mccoy}:

\begin{theorem}
	If $R$ is a commutative ring and $f(x) \in R[x]$ is a zero divisor, then there exists $c \in R$ such that $c f(x) = 0$.
\end{theorem}

 In other words, if $f(x) = a_n x^n + \dots + a_0 \in \zd(\Acal)$, then $c f(x) = c a_n x^n + \dots + c a_0 = 0 \implies c a_n = 0, \dots c a_0 = 0$ by the theorem, which means that $a_n, \dots , a_0$ must be zero divisors in $R$. Note however that this does not mean that if the coefficients of a polynomial in $R[x]$ are all zero divisors that the polynomial will necessarily be a zero divisor. For example, consider $e_1x + e_2 \in \mathbb{R}^2[x]$, which is not a zero divisor in $\mathbb{R}^2$.

Thus, although useful, McCoy's Theorem does not allow us to characterize the zero divisors in polynomial rings built over algebras. However, using the annihilators of the zero divisors of our algebra, we can completely characterize the zero divisors in $\mathcal{A}[z]$ with the following theorem:

\begin{theorem}
	\label{poly zero-divisor char}
	A polynomial $f(z) = a_n z^n + \dots + a_0 \in \mathcal{A}[z]$, where $\mathcal{A}$ is commutative, is a zero divisor if and only if $Ann(a_n) \cap \dots \cap Ann(a_0) \neq \{ 0 \} $
	
	\begin{proof}
		If $f(z)$ is a zero divisor, then by McCoy's theorem, there exists $\eta \in \mathcal{A}$ such that $\eta f(z) = \eta a_n z^n + \dots + \eta a_0 = 0$, which implies that $\eta a_k = 0$ for all $k \in \{ 0, \dots, n \}$. In other words, $\eta$ is in the annihilator of each $a_k$, so $Ann(a_n) \cap \dots \cap Ann(a_0) \neq \{ 0 \}$.
		
		Conversely, suppose $Ann(a_n) \cap \dots \cap Ann(a_0) \neq \{ 0 \}$, then let $\eta \in  Ann(a_n) \cap \dots \cap Ann(a_0) \neq \{ 0 \}$, $\eta \neq 0$. Since $\eta \in Ann(a_k)$ for all $k \in \{ 0, \dots, n \}$, $\eta a_k = 0 \; \forall k \implies \eta a_n z^n + \dots \eta a_0 = \eta f(z) = 0 \implies f(z)$ is a zero divisor.
	\end{proof}
	
\end{theorem}

 This shows, for example, that the polynomial $e_1 x + e_2 \in \mathbb{R}^2[x]$ is in fact not a zero divisor, as we claimed earlier, since $Ann(e_1) \cap Ann(e_2) = \{ 0 \}$. Thus, assuming we have determined the structure of the zero divisors and annihilators of an algebra $\mathcal{A}$, Theorem  
\ref{poly zero-divisor char} gives a relatively quick way to assess whether or not a polynomial in $\mathcal{A}$ is a zero divisor.

Furthermore, if a ring satisfies the Armendariz condition \cite{CAMILLO2008599}, which we define below, we may also characterize the annihilators of elements in $R[x]$:

\begin{definition}
	We say a ring $R$ is Armendariz if for all $p(x), q(x) \in zd(R[x])$, where $p(x) = c_n x^n + \dots + c_1 x + c_0$ and $q(x) = a_m x^m + \dots + a_1 x + a_0$, then $p(x) q(x) = 0$ if and only if $a_i c_j = 0$ for all $i,j$.
\end{definition}

\begin{theorem}
	Let $R$ be an Armendariz ring, and $f(x) \in R[x]$ be a zero divisor, then $g(x) \in \mathrm{Ann}(f(x))$ if and only if $$a_k \in \bigcap_{i=0}^{m} \mathrm{Ann}(b_i)$$ for all $k \in \{ 1, \dots, n \}$ where $g(x) = a_n x^n + \dots + a_0$ and $f(x) = b_m x^m + \dots + b_0$.
	
	\begin{proof}
		If $R$ is Armendariz, then we have $f(x)g(x) = 0$ if and only if $a_k b_j = 0$ for all $k,j$, but by definition we have $f(x)g(x) = 0$ if and only if $g(x) \in \mathrm{Ann}(f(x))$, and $a_k b_j = 0$ if and only if for each $k$ we have $a_k \in \mathrm{Ann}(b_j)$ for all $j$ which again is true if and only if: $$a_k \in \bigcap_{i=0}^{m} \mathrm{Ann}(b_i)$$ for all $k$. 
	\end{proof}
	
\end{theorem}

 An important theorem with regard to Armendariz rings is that all reduced rings (that is, rings without nilpotent elements, and hence in our context semisimple algebras) are Armendariz rings. Furthermore, the class of nil numbers $\mathbf{\Gamma}_n$ are also Armendariz. This first theorem is proven in the introduction to Armendariz rings \cite{Kim00}, and the second claim is a simple corollary of the result proven by Anderson \cite{Anderson98} that $R[x]/\langle x^n \rangle$ is an Armendariz ring if and only if $R$ is a reduced ring for $n \geq 2$, and hence, we know that $\mathbf{\Gamma}_n$ is Armendariz for all $n \geq 2$.

We may also provide a somewhat nicer classification of the zero divisors in a polynomial algebra over $\mathcal{A}$ provided the algebra satisfy additional constraints. We call algebras with this stronger characterization of zero divisors in the polynomial ring \textit{nilfactorable}, which we define below:

\begin{definition}
	If a ring $R$ has the property that for all zero divisors $f(x) \in R[x]$, $f(x) = \epsilon g(x)$ for some zero divisor $\epsilon \in R$, and some non-zero divisor polynomial $g(x) \in R[x]$, then we say that $\mathcal{A}$ is a \textit{nilfactorable} algebra.
\end{definition}

 This is in fact a rather strong property, however we will show that it holds both in $\mathcal{H}$, and in $\mathbf{\Gamma}_n$ for all $n$ by establishing a condition sufficient to guarantee than an algebra is nilfactorable.

\begin{lemma}
	\label{nilfactor lemma}
	If a commutative algebra $\mathcal{A}$ has the property that for all zero divisors $\eta \in \mathcal{A}$, $Ann(\eta) = \xi X$ for some other zero divisor $\xi \in \mathcal{A}$, and some set $X \subseteq \mathcal{A}$, then $\mathcal{A}$ is a \textit{nilfactorable} algebra.
\end{lemma}

 Thus, $\mathbf{\Gamma}_n$ is nilfactorable for all $n$, since as the reader should confirm, for all $\zeta \in \zd(\mathbf{\Gamma}_n)$, $\Ann(\zeta)$ factors as $\epsilon X$ for some $X \subseteq \mathbf{\Gamma}_n$.\footnote{We invite the reader to attempt to see this for themselves, as the structure of the zero divisors in $\mathbf{\Gamma}_n$ is particularly simple. Alternatively, we briefly discuss a technique for characterizing the zero divisors in large class of unital nil algebras, including $\mathbf{\Gamma}_n$, in Section \ref{sec:nil_poset} .} Showing that $\mathcal{H}$ is nilfactorable involves a very simple corollary of Lemma \ref{nilfactor lemma}, which we prove below:

\begin{definition}
	Recall that the annihilator of an element $a \in \mathcal{A}$ of a commutative algebra, denoted $\mathrm{Ann}(a)$ is the set of all elements $x \in \mathcal{A}$ such that $xa = 0$. We claim that this set forms a vector space over $\mathbb{R}$, and call the dimension of this vector space the nildegree of the element $a$. We denote this by $\mathrm{Nil}(a) = \dim(\mathrm{Ann}(a))$.
\end{definition}

\begin{corollary}
	If a finite dimensional commutative semisimple algebra $\mathcal{A}$ has the property that every zero divisor $\eta \in \mathcal{A}$ has nildegree 1, then every polynomial in $\mathcal{A}[x]$ may be factored as the product of a zero divisor in $\mathcal{A}$ and another polynomial in $\mathcal{A}[x]$.
	
	\begin{proof}
		If every zero divisor $a \in \mathcal{A}$ has nildegree 1, $Ann(a)$ is one dimensional, so $Ann(a) = \mathrm{span}\{e\} = e \mathbb{R}$ for some $e \in \mathcal{A}$, and hence, by Lemma \ref{nilfactor lemma}, $\mathcal{A}$ is nilfactorable.
	\end{proof}
	
\end{corollary}

 Importantly, since every finite dimensional algebra is Noetherian, if $f(z) \in \mathcal{A}[z]$ is a zero divisor, and $\mathcal{A}$ is a nilfactorable algebra, then in the decomposition $f(z) = \epsilon_1 g_1(z)$, if $g_1(z)$ is again a zero divisor, we may continue the decomposition $f(z) = \epsilon_1 \epsilon_2 g_2(z)$, and so on if $g_2(z)$ is still a zero divisor. This process gives us the increasing chain of ideals $\langle f(z) \rangle \leq \langle g_1(z) \rangle \leq \dots$, which by the ascending chain condition on ideals must be a finite chain, and hence this process terminates with $f(z) = \epsilon_1\epsilon_2 \dots \epsilon_n g_n(z)$ where $g_n(z)$ is not a zero divisor, therefore:

\begin{proposition}
	\label{prop:nilfact_decomp}
	If $\mathcal{A}$ is a finite dimensional nilfactorable algebra, and $f(z) \in \mathcal{A}[z]$ is a zero divisors, then $f(z) = \epsilon g(z)$ for some $\epsilon \in \mathrm{zd}(\mathcal{A})$ and some $g(z) \in \mathcal{A}[z]$ which is not a zero divisor.
\end{proposition}

 This then allows us to characterize the zero divisors of polynomials over nilfactorable algebras as follows:

\begin{theorem}
	If $\mathcal{A}$ is a finite dimensional nilfactorable algebra, and $f(z) \in \mathcal{A}[z]$ is a zero divisor, with $f(z) = \epsilon g(z)$ the decomposition given in Proposition \ref{prop:nilfact_decomp}, then $\mathrm{Ann}(f(z)) = \mathrm{Ann}(\epsilon) \mathcal{A}[z]$.
	
	\begin{proof}
		Let $\xi \in \mathrm{Ann}(\epsilon)$, $h(z) \in \mathcal{A}[z]$, then $\xi h(z) f(z) = \xi h(z) (\epsilon g(z)) = \xi \epsilon h(z) g(z) = 0$, so $\mathrm{Ann}(\epsilon) \mathcal{A}[z] \leq \mathrm{Ann}(f(z))$.
		
		Conversely, suppose $h(z) \in \mathcal{A}[z]$ and $f(z) h(z) = 0$, then $h(z)$ is a zero divisor, so by Proposition \ref{prop:nilfact_decomp} there exists $\xi \in \mathrm{zd}(\mathcal{A})$ and $k(z) \in \mathcal{A}[z]$ which is not a zero divisor such that $h(z) = \xi k(z)$. Thus, $f(z) h(z) = \epsilon g(z) \xi k(z) = \epsilon \xi g(z) k(z) = 0$ implies $\epsilon \xi = 0$, since neither $g(z)$ nor $k(z)$ are zero divisors, and hence, $\xi \in \mathrm{Ann}(\epsilon)$. Thus, we have shown that $\mathrm{Ann}(f(z)) \leq \mathrm{Ann}(\epsilon) \mathcal{A}[z]$. Therefore, $\mathrm{Ann}(\epsilon) \mathcal{A}[z] = \mathrm{Ann}(f(z))$.
	\end{proof}
\end{theorem}
	
\subsection{The Nil Poset}
\label{sec:nil_poset}

In the last section, we completely characterized the zero divisors for a semisimple algebra. It turns out, as seems to be a common theme in all our work, that the non-semisimple case is more difficult. While some of the techniques of the preceding section will still be useful in this case, studying the zero divisors in more general unital nil algebras (i.e. non-nilfactorable) contexts will likely require more machinery.

In this section, we present one technique which can be used to study the zero divisors of a large class of nil algebras. 

\begin{definition}
	Let $\mathcal{A}$ be an algebra. We say a basis $\beta = \{ v_1, \dots , v_n \}$ of $\mathcal{A}$ is multiplicative if for all $v_i,v_j \in \beta$ we have $v_i \star v_j = c v_k$ for some $c \in \mathbb{R}$ and $v_k \in \beta$. If an algebra $\mathcal{A}$ admits a multiplicative basis, then we say $\mathcal{A}$ is a multiplicative algebra.
\end{definition}

 Given such an algebra which is also unital nil, there is a natural ordering which we may define on the basis:

\begin{definition}
		\label{def:nil_poset}
		Given a unital nil algebra $\mathcal{A}$ with multiplicative basis $\{v_1, v_2, \dots , v_n\}$, where we take $v_1 = 1$ without loss of generality. Define the set $\mathcal{N}_\mathcal{A} = \{0, 1, v_1, \dots, v_n \}$. Then, define an ordering on $\mathcal{N}_\mathcal{A}$ by setting $v_i \preceq v_j$ if and only if there exists $v_k$ and $c \in \mathbb{R}$ such that $v_i \star v_k = c v_j$. As is common in the order theory literature, we will also use $v_i \prec v_j$ as shorthand for $v_i \preceq v_j \wedge v_i \neq v_j$.
\end{definition}

 From the properties of multiplicative nil bases, and the definition of the nil poset, we leave it as an easy exercise to show that:
	
\begin{proposition}
	Given a unital nil algebra $\mathcal{A}$ with a multiplicative basis, $\mathcal{N}_\mathcal{A}$ is a poset. Hence, we call $\mathcal{N}_\mathcal{A}$ the \textit{nil poset} of $\mathcal{A}$ with respect to the basis $\beta$.
\end{proposition} 

 While this structure will also play a role in our study of logarithms, which is where this structure was first conceived, we note that it also has a general role in characterizing the annihilators of an algebra. For instance, by following all possible paths from basis elements to the zero node on the Hasse diagram of a nil poset, we can deduce the annihilators of basis elements in the algebra. To illustrate this, we give some examples of Hasse diagrams of nil posets of some simple nil algebras in the figure below:

\begin{figure}[h!]
	\centering
	\begin{tikzpicture}
	\node (A) {$1$};
	
	\node (S) [below of= A, node distance = 1cm] {};
	
	\node (B) [above of = A, node distance = 1cm] {$\epsilon$} ;
	\node (C) [above of = B, node distance = 1cm] {$0$};
	\draw (A) -- (B) -- (C);
	\end{tikzpicture} \space\space\space\space\space\space\space\space\space\space
	\begin{tikzpicture}
	\node (A) {$1$};
	\node (B) [left of=A,above of=A,node distance=1cm] {$\epsilon$};
	\node (C) [right of=A,above of=A,node distance=1cm] {$\gamma$};
	
	\node (D) [above of=B,node distance = 1.6cm] {};
	\node (D')[left of=D, node distance = .6cm] {$\epsilon^2$};
	
	\node (E) [above of=C,node distance = 1.6cm] {};
	\node (E')[right of=E,node distance = .6cm] {$\gamma^2$};
	
	\node (F) [above of=E,left of=E,node distance = 1.1cm] {};
	\node (F') [above of = F, node distance = .35 cm] {$0$};
	\node (G) [above of=A,node distance = 1.75cm] {$\epsilon \gamma$};
	
	\node (H) [above of=E,left of=E, node distance = .5 cm] {};
	\node (H') [left of = H, node distance = .1 cm] {$\epsilon\gamma^2$};
	\node (I) [above of =D, right of = D, node distance = .5 cm] {};
	\node (I') [right of = I, node distance = .1 cm] {$\epsilon^2 \gamma$};
	
	\draw (A) -- (B) -- (D') -- (F');
	\draw (A) -- (C) -- (E') -- (F');
	\draw (B) -- (G) -- (C);
	\draw (G) -- (H') -- (F') -- (I') -- (G);
	\draw (E') -- (H');
	\draw (I') -- (D');
	\end{tikzpicture} \space\space\space\space\space\space
	\begin{tikzpicture}
	\node (A) {$0$};
	\node (B) [node distance = 1cm, right of = A, below of = A] {$\gamma$};
	\node (C) [node distance = 1cm, left of =A, below of =A] {$\epsilon$};
	\node (D) [node distance = 1.9cm, below of= A] {$1$};
	\node (S) [below of=D, node distance = 1cm] {};
	\draw (A) -- (B) -- (D) -- (C) -- (A);
	\end{tikzpicture}
	
	\caption{Hasse diagrams of the nil posets of $\Gamma_3, \Gamma_3 \otimes \Gamma_3$, and $\Xi_3$, respectively}
\end{figure}
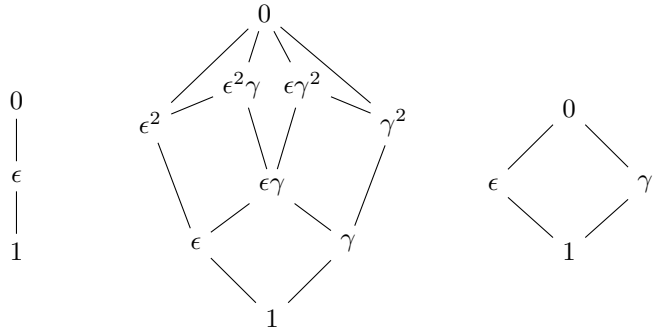

 Notice that each non-zero node in the Hasse diagram of $\mathcal{N}_{\Gamma_3 \otimes \Gamma_3}$ has either one or two elements that cover it. If $z$ is the node in question, these correspond to the elements $\epsilon z$ and $\gamma z$ (if there is only one element that covers the node, as is the case for $\epsilon^2\gamma$ and $\epsilon \gamma^2$, both of these elements coincide, but there are still two distinct left and right paths for $\epsilon z$ and $\gamma z$ which correspond to the single line drawn in the Hasse diagram). Thus, starting from a node $z$ in the Hasse diagram and following a sequence of covering relations left and right corresponds to an element $\epsilon^n \gamma^m $ which annihilates $z$, where $n$ is the number of left steps taken in the path, and $m$ is the number of right steps taken. For example, from this Hasse diagram, we can read $\mathrm{Ann}(\epsilon\gamma) = \mathrm{span}\{ \epsilon^2, \gamma^2, \epsilon \gamma \}$.

The astute reader may note that all of the examples given above are in fact \textit{bounded lattices}. However, this is not true in general. We need simply to consider the algebra $$\mathbb{R}[\epsilon,\gamma,\delta,\eta,\zeta,\xi]/\langle \epsilon^2, \gamma^2, \delta^2,\eta^2, \epsilon\delta - \zeta, \epsilon\eta - \xi, \gamma\delta - \xi, \gamma\eta - \zeta \rangle$$ for a counterexample to this claim. As the reader can easily verify, this algebra presentation is non-degenerate, and the elements $\epsilon,\gamma$ have two distinct minimal upper bounds, $\xi$ and $\zeta$.

Beyond the basic visual technique whereby we may use this order to quickly read off the annihilators of an algebra, this structure leads to a number of interesting open questions. How do various order theoretic properties of the nil poset, such as being a lattice, modularity, and distributivity relate to properties of the algebra? Can any poset be represented as a nil poset, or are there certain posets that cannot be? If some posets cannot be represented as nil posets, is there a simply criterion that characterizes this?

While these are interesting questions, and their relationship to our broader research goals are not at the moment clear, we leave these questions, and other inquiries related to the nil poset open to future researchers interested in our program. 

\pagebreak
\bibliography{sources}

\end{document}